\documentclass[11pt]{article}
\usepackage{amsmath}
\usepackage{amssymb}
\usepackage{txfonts}
\usepackage{lscape}
\usepackage{algpseudocode}
\usepackage{MnSymbol}
\title{A note on banded linear systems}
\author{
D. Barrios Rolan\'{\i}a, \, J.C. Garc\'{\i}a-Ardila     \\
Depto. Matem\'atica Aplicada a la Ingenier\'{\i}a Industrial\\
Universidad Polit\'ecnica de Madrid 
}
\date{}
\begin{document}

\maketitle

 \begin{abstract} 
In \cite{Barrios-Manrique2015}  a new factorization for infinite Hessenberg banded matrices was introduced. In this note we prove that this kind of factorization can also be used for finite matrices. In addition, a new method for solving banded linear systems is provided.
 \end{abstract}

\section{Introduction}

Banded linear systems constitute a relevant kind of linear systems in scientific computing due to its applications in many areas of science and engineering. These systems arise in the study of $p$-orthogonal polynomials and other fields of Approximation Theory, as well as the discretization and linearization of differential equations \cite{Libro}, \cite{Saugata}. In particular, tridiagonal linear systems are asociated to cubic splines, quadrature formulas an other subjects where the zeros of a sequence of orthogonal polynomials have to be located \cite{Lotar}.

An extensive class of direct methods for solving a linear system
\begin{equation}
A_NX=b
\label{finito}
\end{equation}
is based on the 
$LU$
triangular decomposition 
\begin{equation}
A_N=L_NU_N
\label{finito2}
\end{equation}
of the coefficient matrix \cite{Duff,Tref97}. It is known that there is no universally best method for solving linear systems. In fact, the choice of one or the other method depends on the problem under consideration, which justifies the construction of new methods in addition to the already known ones. In this sense we emphasize that, if  $A_N$ is a Hessenberg matrix, not necessarily banded, there are several sophisticated methods to deal with \eqref{finito} (see for example \cite{Jandron,Ru11,Stone 73}).  

It is well known that, when $A_N$ is a finite Hessenberg banded matrix of order $N$, 
$$
A_N = \left( \begin{array}{cccccc}
 a_{0,0} & a_{0,1}& 0&\cdots  &\cdots & 0    \\
 a_{1,0} & a_{1,1} & a_{1,2} && &\vdots \\
 \vdots      & \vdots    &  \ddots   & \ddots  &&\vdots     \\
 a_{p,0} &a_{p,1} & \cdots & a_{p,p} &a_{p,p+1}& \vdots \\
0& a_{p+1,1} & &\ddots& \ddots &\vdots  \\
\vdots & 0 &&\ddots &\ddots &\vdots   \\
\vdots  &\ddots &\ddots &\ddots &\ddots &0\\
\vdots &\ddots&0 &a_{N-2,N-p-1}&\cdots & a_{N-2,N-1}\\
0 & \cdots & 0 & a_{N-1,N-p-1}& \cdots & a_{N-1,N-1}
\end{array} \right)\,,
$$
with 
$
a_{i+p,i}\neq 0,\,i=0,1,\ldots,N-p-1, 
$
and the factorization \eqref{finito2} can be obtained, then 
$
L_N
$
and 
$U_N$ 
are triangular banded matrices of the same order $N$. We assume $N>>p$.  Further, the diagonal entries of $L_N$ can be assumed equal to 1, being
\begin{equation} 	
\label{L}
L_N=
\left(
\begin{array}{ccccccccccc}
1        &      &             &\\
l_{1,0}&1     &             & \\
\vdots & \vdots  &\ddots&\\
\vdots & \vdots & \ddots  &\ddots&\\
l_{p,0}&l_{p,1}  &\cdots&l_{p,p-1}&1\\
0      &l_{p+1,1}&   \cdots  &  \cdots& l_{p+1,p}& 1\\
\vdots      & \ddots    &  \ddots   && \ddots    &\ddots    \\
0      & \cdots    &0 &  l_{N-1,N-p-1} &\cdots   & l_{N-1,N-2}     & 1 \\
\end{array}
\right),
\end{equation}	
with 
$l_{p+i,i}\ne 0$ for $i=0,1,\ldots,N-p-1$. In this case $U_N$ is a bidiagonal upper triangular matrix,
\begin{equation} 	
\label{U}
U_N=
\left(
\begin{array}{lllllcccc}
u_{1}&1     &             & \\
& u_{2} &\ddots&\\
&& \ddots  &\ddots&\\
&   && u_{N-1}    &1    \\
&&&& u_{N}\\
\end{array}
\right).
\end{equation}
Under the above conditions, matrices $L_N$ and $U_N$ are uniquely determined.

On the other hand, the Darboux factorization for an infinite lower triangular $(p+1)$-banded matrix $L$ was introduced and analyzed in \cite{Barrios2015}. Assuming
\begin{equation} 	
\label{Linfinita}
L=
\left(
\begin{array}{ccccccccccc}
1        &      &             &\\
l_{1,0}&1     &             & \\
\vdots & \vdots  &\ddots&\\
\vdots & \vdots & \ddots  &\ddots&\\
l_{p,0}&l_{p,1}  &\cdots&l_{p,p-1}&1\\
0      &l_{p+1,1}&   \cdots  &  \cdots& l_{p+1,p}& 1\\
\vdots      & \ddots    &  \ddots   && \ddots    &\ddots   
\end{array}
\right)
\end{equation}	
and
$l_{p+i,i}\ne 0$ for $i=0,1,\ldots$, the existence of $p$ bi-diagonal infinite matrices $L^{(i)}$, $i=1,2,\ldots,p,$
$$
	L^{(i)}=
\left(
\begin{array}{ccccl}
	1&&& \\
	\gamma_{i+1}&1&&\\
	&\gamma_{p+i+2}&1&\\
	&   &   \gamma_{2p+i+3}      &\ddots\\
& & & \ddots
	\end{array}
\right), 
\quad \gamma_{jp+i+j+1}\ne 0,\quad j=0,1,\ldots,
$$
verifying 
\begin{equation}\label{factorization}
L=L^{(1)}L^{(2)}\cdots L^{(p)}	
\end{equation}
was proved. This decomposition is no unique since it depends on the choice of the set of entries
\begin{equation}\label{e.arbitrary}
\begin{array}{clllll}
\gamma_2&\cdots &\gamma_{p-1}&\gamma_p\\
\gamma_{p+3}& \cdots &\gamma_{2p}\\
\vdots &\udots &\\
\gamma_{p(p-1)}\\
\end{array}
\end{equation}
(see Table \ref{tabla}). In this note we show that this factorization can be used in the case of finite matrices. As a consequence, a new method for solving a linear system \eqref{finito} is provided.

\section{Darboux factorization for finite matrices}

For an infinite lower banded matrix $A$, we assume $A=LU$ where $L$ is given as in \eqref{Linfinita} and $U$ is an upper triangular bidiagonal matrix, this is,
$$
U=
\left(
\begin{array}{ccccl}
	\gamma_{1}&1&&\\
	&\gamma_{p+2}&1&\\
	&   &   \gamma_{2p+3}      &\ddots\\
& & & \ddots
	\end{array}
\right).
$$
Following \cite{Barrios2015}, we assume $L$ decomposed as in \eqref{factorization}. Each row of Table \ref{tabla} represents the corresponding rows of $U$ and the set of factors of \eqref{factorization}. 

We consider the secondary diagonal in Table \ref{tabla}, this is,
$$
\gamma_{p+1},\,\gamma_{2p+1},\ldots, \gamma_{p^2+1},\,\gamma_{(p+1)p+1}.
$$
At the top of this secondary diagonal we see some entries framed, which are the starting data \eqref{e.arbitrary}. Furthermore, the entries of $U$, at the first column of Table \ref{tabla}, are well known from the $LU$ factorization. Our main aim in this section is to show that each one of the rest of entries can be determine from the previous rows.

\begin{landscape}
\begin{table}
\caption{Factors of $L$}
\label{tabla}
\begin{tabular}{cccccccccccccc}
 $U$ &   $L^{(1)}$ & $L^{(2)}$  & $\cdots$ &  $L^{(p-s)}$ &$\cdots $&$L^{(p-2)}$& $L^{(p-1)}$& $L^{(p)}$& \\\\
\hline\\
$ \gamma_{1}$ & $\boxed{\gamma_{2}}$ &$\boxed{\gamma_3}$ & $\cdots$ &  $\boxed{\gamma_{p-s+1}}$ & $\cdots$ & 
 $\boxed{\gamma_{p-1}}$ &  $\boxed{\gamma_{p}}$ & $\gamma_{p+1}$ \\\\
 $ \gamma_{p+2}$ & $\boxed{\gamma_{p+3}}$ &$\boxed{\gamma_{p+4}}$ & $\cdots$ &  $\boxed{\gamma_{2p-s+2}}$ & $\cdots$ & 
 $\boxed{\gamma_{2p}}$ &  $\gamma_{2p+1}$ & $\gamma_{2p+2}$\\\\
$ \gamma_{2p+3}$ & $\boxed{\gamma_{2p+4}}$ &$\boxed{\gamma_{2p+5}}$ & $\cdots$ &  $\boxed{\gamma_{3p-s+3}}$ & $\cdots$ & 
 $\gamma_{3p+1}$ &  $\gamma_{3p+2}$ & $\gamma_{3p+3}$\\\\
$\vdots$ & $\vdots$ &$\vdots$&  & $\vdots$ & &$\vdots$  &$\vdots$&$\vdots$  \\\\
$ \gamma_{(s-1)p+s}$ & $\boxed{\gamma_{(s-1)p+s+1}}$ &$\boxed{\gamma_{(s-1)p+s+2}}$ & $\cdots$ &  $\boxed{\gamma_{sp}}$ & $\cdots$ & 
 $\gamma_{sp+s-2}$ &  $\gamma_{sp+s-1}$ & $\gamma_{sp+s}$\\\\
$ \gamma_{sp+s+1}$ & $\boxed{\gamma_{sp+s+2}}$ &$\boxed{\gamma_{sp+s+3}}$ & $\cdots$ &  $\gamma_{(s+1)p+1}$ & $\cdots$ & 
 $\gamma_{(s+1)p+s-1}$ &  $\gamma_{(s+1)p+s}$ & $\gamma_{(s+1)p+s+1}$\\\\
$\vdots$ & $\vdots$ &$\vdots$&  & $\vdots$ & &$\vdots$  &$\vdots$&$\vdots$  \\\\
 $ \gamma_{(p-2)p-2}$ & $\boxed{\gamma_{(p-2)p-1}}$ &$\boxed{\gamma_{(p-2)p}}$ & $\cdots$ &  $\gamma_{(p-1)p-s-2}$ & $\cdots$ & 
 $\gamma_{(p-2)p+p-4}$ &  $\gamma_{(p-2)p+p-3}$ & $\gamma_{(p-2)p+p-2}$\\\\
$ \gamma_{(p-1)p-1}$ & $\boxed{\gamma_{(p-1)p}}$ &$\gamma_{(p-1)p+1}$ & $\cdots$ &  $\gamma_{p^2-s-1}$ & $\cdots$ & 
 $\gamma_{(p-1)p+p-3}$ &  $\gamma_{(p-1)p+p-2}$ & $\gamma_{(p-1)p+p-1}$\\\\
$ \gamma_{p^2}$ & $\gamma_{p^2+1}$ &$\gamma_{p^2+2}$ & $\cdots$ &  $\gamma_{(p+1)p-s}$ & $\cdots$ & 
 $\gamma_{p^2+p-2}$ &  $\gamma_{p^2+p-1}$ & $\gamma_{p^2+p}$\\\\
$\vdots$ & $\vdots$ &$\vdots$&  & $\vdots$ & &$\vdots$  &$\vdots$&$\vdots$  
\end{tabular}
\end{table}
\end{landscape}

We call $s$-th {\it secondary diagonal}, $s=1,2,\ldots,$ the set given by the entries 
$
\gamma_{sp+s},\,\gamma_{(s+1)p+s},\ldots, \gamma_{(s+p-1)p+s},\,\gamma_{(s+p)p+s}
$
in Table \ref{tabla}. In particular, for $s=1$ we have the previously called secondary diagonal. In the following, for each fixed $i\in \mathbb{N}$ we show that the $i$-th secondary diagonal is determined in terms of the previous $s$-th secondary diagonals, $s=1,2,\ldots, i-1$, and the starting data \eqref{e.arbitrary}. Even more, we will see that each entry of this $i$-th secondary diagonal in the $k$-th row is obtained exclusively in terms of such entries that are in the rows $1,2,\ldots, k$.

From \cite[(35)]{Barrios2015} we have 
\begin{eqnarray*}
\delta^{(i)}_{k}\gamma_{(k+i+1)p+i}& = &a_{k+i+1,i-1} - \sum_{
\widetilde{E}^{(0)}_{k+2}} \gamma_{(i-2)p+i+i_1-1}\gamma_{(i-1)p+i+i_2-1}\cdots \gamma_{(k+i)p+i+i_{k+3}-1},\nonumber\\
& & k=-1,0,\ldots , p-2\,,
\label{(a)}
\end{eqnarray*}
where
\begin{equation}
\delta^{(i)}_{k}=\gamma_{(i-1)p+i}\gamma_{ip+i}\cdots \gamma_{(k+i)p+i}
\label{20}
\end{equation}
and
\begin{equation}\label{3434}
\widetilde{E}^{(0)}_{k+2}=\{(i_1,\ldots ,i_{k+3}):k+3\leq i_{k+3}\leq\cdots \leq i_1\leq p+1\,,i_{k+3}<p+1\}.
\end{equation}
For each $k=-1,0,\ldots,p-2$, the entry $\gamma_{(k+i+1)p+i}$ is in the $(i+k+1)$-th row and $i$-th secondary diagonal. Since \eqref{(a)} we can express this entry in terms of 
$\delta^{(i)}_{k}$ and 
\begin{equation}
\label{**}
\gamma_{(i-2)p+i+i_1-1},\,\gamma_{(i-1)p+i+i_2-1},\ldots ,\,\gamma_{(k+i)p+i+i_{k+3}-1}
\end{equation}
when 
$
(i_1,\ldots ,i_{k+3})
\in
\widetilde{E}^{(0)}_{k+2}.
$

Firstly, from \eqref{20} we see that
$\delta^{(i)}_{k}$
is computed from the entries of the same $i$-th secondary diagonal that are in the rows 
$i,i+1,\ldots , i+k$. 

Secondly, we analyze the entries \eqref{**}, this is,
\begin{equation}
\label{***}
\gamma_{(r+i)p+i+i_{r+3}-1},\quad r=-2,-1,\ldots, k.
\end{equation}
If $r\leq k-1$ then, taking into account \eqref{3434},
$$
(r+i)p+(i+1)\leq (r+i)p+i+i_{r+3}-1\leq (r+i+1)p+i.
$$ 
Hence 
$
\gamma_{(r+i)p+i+i_{r+3}-1}
$
is in some row of Table \ref{tabla} before the $(r+i-1)$-th row. Moreover, when 
$
\gamma_{(r+i-1)p+i+i_{r+2}-1}
$
is in the $j$-th column then 
$
\gamma_{(r+i)p+i+i_{r+3}-1}
$
is in the $(j-1)$-th column or some previous column of the following row. Therefore, if 
$
\gamma_{(r+i-1)p+i+i_{r+2}-1}
$
is at the top of the $i$-th secondary diagonal, the same is true for 
$
\gamma_{(r+i)p+i+i_{r+3}-1}.
$
Finally, for $r=k$ in \eqref{***} the situation is similar but now
$
\gamma_{(r+i)p+i+i_{r+3}-1}=\gamma_{(k+i)p+i+i_{k+3}-1}
$
is in the $(k+i-1)$-th row and not in the $i$-th secondary diagonal, because 
$
(r+i)p+i+i_{k+3}-1< (r+i+1)p+i.
$
(Just, the entry of this row in the $i$-th secondary diagonal is that we want to compute.)

In summary, each entry in the $i$-th secondary diagonal of Table \ref{tabla} is obtained with the entries of the previous rows that are at the top of the $i$-th secondary diagonal. Translating this reasoning to matrices 
$
L^{(1)},\ldots , L^{(p)}, U,
$
we deduce that the entry in the row $i$ of $L^{(s)},\,s=1,\ldots, p$, is obtained using only the rows $1,2,\ldots , i$ of 
$
L^{(1)},\ldots , L^{(s)}, U.
$
As a consequence, 
$
\left(L^{(1)}\cdots  L^{(p)} U\right)_n=L_n^{(1)}\cdots  L_n^{(p)} U_n\,,\, n\in \mathbb{N}.
$
In particular, 
$$
\left(L^{(1)}\cdots  L^{(p)} U\right)_N=L_N^{(1)}\cdots  L_N^{(p)} U_N\,.
$$
From this and the well-known fact that 
$
\left( LU\right)_N=L_NU_N
$
we obtain 
\begin{equation}
\left(L^{(1)}\cdots  L^{(p)} \right)_N=L_N^{(1)}\cdots  L_N^{(p)} \,.
\label{otroN}
\end{equation}

\section{Darboux factorization and banded systems}

As a consequence of \eqref{otroN}, it is possible to use the Darboux factorization for finite matrices. In other words, if there exists the $LU$ factorization for the coefficients matrix $A_N$ in the system \eqref{finito} then  we have
\begin{equation}\label{ALU}
A_N=L_N^{(1)}\cdots  L_N^{(p)} U_N
\end{equation}
and we can define
$$
X^{(i)}=
\left\{
\begin{array}{lll}
L^{(i+1)}_N\cdots  L^{(p)}_N U_NX& , & i=1,\ldots ,p-1,\\\\
U_NX & , & i=p.
\end{array}
\right.
$$
Thereby, \eqref{finito} is reduced to solve iteratively the following $p+1$ tridiagonal systems,
\begin{equation}\label{otro14}
\left\{
\begin{array}{lll}
L^{(1)}_N X^{(1)}=b\\\\
L^{(k)}_N X^{(k)}=X^{(k-1)}& , & k=2,\ldots ,p\\\\
U_NX=X^{(p)}. 
\end{array}
\right.
\end{equation}
A remarkable advantage of the proposed method \eqref{otro14} is its low computational complexity. In fact, this is an extension of the method usually used to solve tridiagonal systems based in the $LU$ factorization of the coefficients matrix. 

We assume \eqref{finito2}, where $L_N$ and $U_N$ are given by \eqref{L} and \eqref{U}, respectively. With the purpose to derive an algorithm for obtaining the decomposition \eqref{otroN}, we write the entries of $L_N$ verifying this decomposition, this is,
\begin{equation*}\label{Lelements} 
l_{m,m-k}=
\displaystyle\sum_{1\leq \sigma_1<\cdots <\sigma_k\leq p}\left(\prod_{j=1}^{k}\gamma_{(m-j)p+\sigma_j+m-j+1}\right),\quad k=1,\ldots ,p.
\end{equation*}
We recall that matrix $U_N$ is known from the $LU$ factorization of $A_N$. Therefore, we only need to determine the entries  $\gamma_{(m-1)p+m+i}$, which are, for $m=1,2,\ldots,N-1,$ in each row of matrices $L^{(i)}_N,\, i=1,2,\ldots, p$.  Hence, 
\begin{eqnarray}\label{otro19}
l_{m,m-k}&=&\sum_{\substack{1\leq \sigma_1<\cdots <\sigma_k\leq p\\\sigma_1\ne p-k+1}}\left(\prod_{j=1}^{k}\gamma_{(m-j)p+\sigma_j+m-j+1}\right)\nonumber\\
&+&\gamma_{m(p+1)-k+1}\prod_{j=2}^k \gamma_{(m-j+1)p-k+m+1},\quad k=1,\ldots,p.
\end{eqnarray}
Besides, because we are assuming $a_{m,m-p}\ne 0$, then all the entries of $L^{(s)},\,s=1,\ldots, p,$ are necessarily  nonzero. This is, $\gamma_{(j-1)p+s+j}\ne0$  for $s=1,\ldots, p,\,j=1,2,\ldots$ Thus, defining
\begin{equation}\label{ganma}
\Gamma_s:= \prod_{j=2}^{p-s+1}\gamma_{(m-j)p+m+s},\quad s=1,\ldots, p,
\end{equation}
we have
$\Gamma_s\neq 0$ 
and, from this and \eqref{otro19}, taking 
$s=p-k+1$ for $s=1,\ldots, p,$ we can write
\begin{equation}\label{otro}
\gamma_{(m-1)p+m+s}=\left(l_{m,m-p+s-1}-\displaystyle\sum_{\substack{1\leq \sigma_1<\cdots <\sigma_{p-s+1}\leq p\\\sigma_1\ne s}}\prod_{j=1}^{p-s+1}\gamma_{(m-j)p+\sigma_j+m-j+1}\right)/\Gamma_s
\end{equation}
Thereby, when all the entries  
$\gamma_{(k-1)p+k+s},\, s=1,2\ldots, p,$ $k=1,2,\ldots, m-1,$ have been computed, 
$
\gamma_{(m-1)p+m+s}
$
can be computed using \eqref{ganma}- \eqref{otro}. 

In this way, from the starting data \eqref{e.arbitrary} we get row by row those of Table \ref{tabla}. If $a \to b$ means that $b$ is obtained from $a$, in a schematic form we write  
\begin{align*}
&\to \gamma_{p+1}\\
&\to \gamma_{2p+1}\to\gamma_{2p+2}\\
&\to \gamma_{3p+1}\to\gamma_{3p+2}\to \gamma
_{3p+3}\\
&\vdots\\
&\to\gamma_{p^2+1}\to\gamma_{p^2+2}\to\cdots\to  \gamma_{p^2+p}\\
&\vdots\\
&\to\gamma_{(N-2)p+N}\to\gamma_{(N-2)p+N+1}\to\cdots\to  \gamma_{(N-1)p+N-1}
\end{align*}

For 
$p>1$
the algorithm can be summarized as follows,

\vspace{.5cm}

\noindent
\begin{algorithmic}[1]
	\Require: $(l_{i,j})$ in \eqref{L}. Specify nonzero values for \eqref{e.arbitrary}.
\For {$m=1, \ldots, p-1$}
	\For{$s=1,\ldots, m$}
	\State{\begin{multline*}
		\gamma_{mp+s}=\left(l_{m,s-1}-\sum_{\substack{1\leq \sigma_1<\cdots <\sigma_p\leq p \\\sigma_1\ne p+s-m}}\prod_{j=1}^{m-s+1}\gamma_{(m-j)p+\sigma_j+m-j+1} \right)\left/\prod_{j=2}^{m-s+1}
		\gamma_{(m-j+1)p+s} \right.
		\end{multline*}
	}
	\EndFor
	\EndFor
	\For {$m=p, p+1,\ldots,N-1,$}
	\For{$s=m-p+1,m-p+2,\ldots, m$}
	\State{\begin{multline*}
		\gamma_{mp+s}=\left(l_{m,s-1}-\sum_{\substack{1\leq \sigma_1<\cdots <\sigma_p\leq p \\\sigma_1\ne p+s-m}}\prod_{j=1}^{m-s+1}\gamma_{(m-j)p+\sigma_j+m-j+1} \right)\left/\prod_{j=2}^{m-s+1}
		\gamma_{(m-j+1)p+s} \right.
		\end{multline*}
	}	
	\EndFor
	\EndFor
		\end{algorithmic}

\vspace{.5cm}

The work involved in the second loop (line 2$\to$line 4) is dominated by 
\begin{multline*}
		\gamma_{mp+s}=\left(l_{m,s-1}-\sum_{\substack{1\leq \sigma_1<\cdots <\sigma_p\leq p \\\sigma_1\ne p+s-m}}\prod_{j=1}^{m-s+1}\gamma_{(m-j)p+\sigma_j+m-j+1} \right)\left/\prod_{j=2}^{m-s+1}
		\gamma_{(m-j+1)p+s} \right.
		\end{multline*}
This computation requires 
$\binom{p}{m-s+1}(m-s)+1$
scalar operations, of which $[\binom{p}{m-s+1}-1](m-s)$ are sums and products of the numerator, 1 subtraction and 1 division, as well as  $(m-s-1)$ products in the denominator. Then the number of operations in lines $1$ to $5$ is
$$
M_1=\sum_{m=1}^{p-1}\sum_{s=1}^m\left(\binom{p}{m-s+1}(m-s)+1\right)=\frac{p(p-1)}{2}+p!\sum_{s=2}^{p-1}\frac{s-1}{s!(p-s-1)!},
$$
where we understand that the last term is 0 when $p=2$. It is easy to see that 
$$
\sum_{s=1}^{p-1}\binom{p}{s-1}s(p-s-1)=p+(p-3)p2^{p-2}.
$$
Hence
$$
M_1=\frac{p}{2}(p+1+(p-3)2^{p-1}).
$$
Doing a similar analysis, lines 7 to 9 require 
$\binom{p}{m-s+1}(m-s)+1$
scalar operations for 
$s=m-p+1, m-p+2,\ldots, m$. Then, in lines 6 to 10 we have
$$
M_2=p(N-p)+\sum_{m=p}^{N-1}\left[\sum_{s=m-p+1}^m\binom{p}{m-s+1}(m-s)\right]=(N-p)(1+(p-2)2^{p-1})
$$
operations. Therefore, the number of operations for the factorization \eqref{otroN} is
$$
M_1+M_2=(1-2^{p-1})\frac{p(p-1)}{2}+(1+(p-2)2^{p-1})N.$$

When
$A_N$
is a Hessenberg matrix, from the above it is immediate to build an analogous algorithm including the matrix U. In other words, it is possible to find the factorization \eqref{ALU} without to know the previous $LU$ factorization of $A_N$. However, if $A_N$ is a $q$-banded matrix, but no a Hessenberg matrix, its $LU$ factorization drives to a banded upper triangular matrix $U_N$ and, making use of the same idea of this section, 
$$
U_N^T={U_N^{(q)}}^T\cdots {U_N^{(1)}}^T.
$$
Then 
$$
A_N=L_N^{(1)}\cdots L_N^{(p)}U_N^{(1)}\cdots U_N^{(q)}
$$
and \eqref{otro14} can be applied in a more general situation.

\end{document}